\documentclass[12pt]{article}
\usepackage{amsmath}
\usepackage{amssymb}
\usepackage{dsfont}
\usepackage{cite}
\newsymbol \blackbox 1004
\newcommand{\eh}{\hfill}\newlength{\sperr}

\newenvironment{proof}{{\settowidth{\sperr}{\bf\rm
Proof}%
\par\addvspace{0.3cm}\noindent\parbox[t]{1.3\sperr}
{\bf\rm P\eh r\eh o\eh o\eh f\eh }%
}}{\nopagebreak\mbox{}
$\blackbox$\par\addvspace{0.3cm}}

\def\nn{\nonumber}

\def\a{\alpha}
\def\b{\beta}

\def\g{\gamma}

\def\Lam{\Lambda}
\def\s{\sigma}
\def\la{\lambda}
\def\om{\omega}

\def\t{\theta}
\def\ze{\zeta}

\def\vp{\varphi}
\def\vt{\vartheta}

\def\wh{\widehat}
\def\wt{\widetilde}
\def\ov{\overline}

\def\BC{{\mathbb C}}
\def\BD{{\mathbb D}}
\def\BR{{\mathbb R}}
\def\BN{{\mathbb N}}

\def\cla{{\mathcal A}}

\def\clc{{\mathcal C}}

\def\clm{{\mathcal M}}
\def\cln{{\mathcal N}}
\def\clu{{\mathcal U}}
\def\clw{{\mathcal W}}

\def\clr{\mathcal{R}}
\def\cls{\mathcal{S}}

\newcommand{\E}{\mathrm{e}}
\newcommand{\I}{\mathrm{i}}
\def\mf{\mathfrak}

\newtheorem{Pa}{Paper}[section]
\newtheorem{Tm}[Pa]{{\bf Theorem}}

\newtheorem{Rk}[Pa]{{\bf Remark}}

\newtheorem{Dn}[Pa]{{\bf Definition}}

\title{Discrete self-adjoint Dirac systems: asymptotic relations, Weyl functions and Toeplitz matrices}

\author{Alexander Sakhnovich}

\date{}
\begin{document}
\maketitle


\begin{abstract}   We consider discrete Dirac systems as an alternative (to the famous
Szeg{\H o} recurrencies and matrix orthogonal polynomials) approach to the study of the
corresponding block Toeplitz matrices. We prove an analog of the Christoffel--Darboux formula
and derive the asymptotic relations for the  analog  of reproducing kernel
(using Weyl--Titchmarsh functions of discrete Dirac systems). We  study also the case of rational Weyl--Titchmarsh functions (and GBDT 
version of the B\"acklund-Darboux transformation
of the trivial discrete Dirac system). We show that block diagonal plus semi-separable Toeplitz matrices appear in this case.
\end{abstract}


\vspace{0.2em}


\section{Introduction}\label{Intro}
\setcounter{equation}{0}

Self-adjoint discrete Dirac systems have been introduced in \cite{FKRS08} and studied further
in \cite{FKRS19, RS18, ALS-JAT, SaSaR} following the case of skew-self-adjoint discrete Dirac
systems in \cite{KS}.  In particular, the paper \cite{ALS-JAT} is dedicated to the interrelations
between self-adjoint discrete Dirac systems and block Toeplitz matrices, which (in the scalar case)
are in many respects similar to the interrelations between the famous Szeg{\H o}  recurrences and Toeplitz matrices.
The corresponding Verblunsky-type theorems are proved in \cite{ALS-JAT}. Our present paper may be considered
as an important development of the work started in \cite{ALS-JAT}. (See also \cite[Chapter 5]{SaSaR} on discrete
Dirac systems.) Some related recent research and references on continuous systems
and convolution and other structured operators 
one can find in \cite{MakP, Polt, SaSaR, SaL2012}

Self-adjoint discrete Dirac system on the semi-axis $0 \leq k<\infty$  is the system of the form 
\begin{align}& \label{I1}
y_{k+1}(\la)=\left(I_{2p}-\frac{\I}{\la}j C_k\right)y_k(\la), \quad j:=\begin{bmatrix}I_p & 0 \\ 0 & -I_p \end{bmatrix},
\end{align}
where $k\in \BN_0=\BN \cup \{0\}$, $\BN$ is the set of positive integers,
 $\la$ is the so called spectral parameter, $y_{k+1}(\la)$ and $y_k(\la)$ are $2p \times 1$ vector functions, $I_p$ is the $p\times p$ identity matrix, 
 $\I$ stands for the  imaginary unit ($\I^2=-1$),
and the $2p \times 2p$ matrices $C_k$ have the following properties:
\begin{align} & \label{I2}
C_k>0, \quad C_k j C_k=j \quad (0 \leq k<\infty).
\end{align}
System \eqref{I1}, \eqref{I2} may be also considered  as a special case of the discrete canonical system.

Here, we introduce  analogs of  Christoffel-Darboux formula  and of reproducing kernels for system \eqref{I1} such that
\eqref{I2} holds. The asymptotics of reproducing kernels  is of essential interest in theory and applications to  random processes
(see, e.g., \cite{BreD, Lub} and references therein). In our paper, we study the asymptotics of the analogs of reproducing kernels
using Weyl--Titchmarsh (Weyl)  functions of system \eqref{I1}. Some of the possible applications are connected with the physical Gaussian models \cite{ALS87, VVol}.

We note that  interesting and related papers on analogs of Weyl functions in the theory of orthogonal polynomials were written, for instance, by L.~Golinskii and P. Nevai \cite{GolN}
and by B.~Simon \cite{SiB}. The case of matrix orthogonal polynomials is of interest as well (see, e.g., \cite{MaMa, DaS, Delsarte, Koe}) but the corresponding analogs of Szeg{\H o} recurrences
are rather complicated whereas discrete Dirac systems, which we consider here, have the same form as in the scalar case. 

Section \ref{Prel} is dedicated to some basic preliminary results in order to make the paper self-contained.
In Section \ref{CD}, we consider  Christoffel-Darboux formula  and asymptotics of the analogs of reproducing kernels.
Finally, in Section \ref{Expl} we study the case of rational Weyl functions and our GBDT version of the B\"acklund-Darboux transformation
of the trivial discrete Dirac system. (See \cite{Ci,  GeT, Gu, KoSaTe, MS, SaSaR, ZM} on B\"acklund-Darboux transformations and related
commutation methods.) We show that block diagonal plus block semiseparable Toeplitz matrices appear in this case.

As usual, $\BN$ stands for the set of positive integers, $\BR$ stands for the real axis, $\BD$ stands for the unit disk $\{z: \,|z|<1\}$,
$\BC$ stands for the complex plane and $L^1(\BR)$ denotes the class of absolutely
integrable functions on $\BR$. The open upper (lower) half-plane is denoted
by $\BC_+$ ($\BC_-$), and $\ov{\BC_+}$ ($\ov{\BC_-}$) stands for the closed upper (lower) half-plane.
The notation $\ov{\la}$ means complex conjugate of $\la$ and $A^*$ means complex conjugate transpose of the matrix $A$.
The inequality $S>0$ for some matrix $S$ means that $S$ is positive definite.
The notation $\Im(\a)$ stands for the imaginary part of matrix $\a$ (i.e.,  $\Im(\a)=\frac{1}{2\I}(\a-\a^*)$).
Below in the text, we write sometimes Dirac system or discrete Dirac system meaning self-adjoint discrete Dirac system
\eqref{I1}, where \eqref{I2} holds. 

\section{Preliminaries}\label{Prel}
\setcounter{equation}{0}
The fundamental solution of the Dirac system \eqref{I1} (where \eqref{I2} holds) is denoted by $\{W_k(\la)\}$,
that is,  $W_{k+1}(\la)=\left(I_{2p}-\frac{\I}{\la}j C_k\right)W_k(\la)$. This solution is normalised by
\begin{align}& \label{P1}
W_0(\la)\equiv I_{2p}.
\end{align}
Let us recall the definition of the Weyl function of  system \eqref{I1} (see, e.g., \cite{ALS-JAT}).
\begin{Dn}\label{DnWD}
 A $p \times p$ matrix function
$\vp(\la)$
holomorphic in the lower complex half-plane $\BC_-$ is called a Weyl function
for Dirac system \eqref{I1}, \eqref{I2}
  if the
inequality
\begin{equation} \label{P2}
\sum_{k=0}^\infty [\I \vp(\lambda)^* \quad I_p
]q(\lambda)^k K
W_k(\lambda)^*C_kW_k(\lambda)K^*\left[\begin{array}{c}
-\I \vp(\lambda) \\ I_p
\end{array}
\right]<\infty ,
\end{equation}
holds for
\begin{align} & \label{P3}
q(\lambda):=|\lambda^2|(|\lambda^2|+1)^{-1}, \quad K:= \frac{1}{\sqrt{2}}\begin{bmatrix}I_p & -I_p \\ I_p & I_p \end{bmatrix}.
\end{align}
\end{Dn}
Weyl functions defined by the inequality \eqref{P2}  belong to Herglotz class, that is, $\Im\big(\vp(\la)\big)\leq 0$ for $\la \in \BC_-$ (see \cite[Corollary 5.5]{FKRS08}).
Formula (6.13) and Theorems 5.2 and 6.5 in \cite{FKRS08} imply that the  Weyl function $\vp$ of Dirac system \eqref{I1}, \eqref{I2} is unique and
admits Taylor representation
\begin{align}& \label{P4}
\I \vp\Big(\I \frac{z+1}{z-1}\Big)=\a_0+\sum_{k=1}^{\infty}s_{-k}z^k,
\end{align}
where (putting $s_0=\a_0+\a_0^*$ and $s_k=s_{-k}^*$) we have
\begin{align}& \label{P5}
S(N)=\{s_{k-i}\}_{i,k=1}^N>0 \quad (1 \leq N <\infty).
\end{align}
Recall that Toeplitz matrices satisfy matrix identities (see a detailed discussion and references in \cite{ALS-JAT}):
\begin{align}& \label{P6}
AS(N)-S(N)A^*= \I \Pi J \Pi^*;  \quad \Pi=\begin{bmatrix}\Phi_1 & \Phi_2\end{bmatrix},  
\end{align}
where $A=A(N), \quad \Pi=\Pi(N), \quad \Phi_1=\Phi_1(N), \quad \Phi_2=\Phi_2(N)$,
\begin{align}&  \ \label{P6+}
A(N)= \left\{ a_{k-i}^{\,} \right\}_{i,k=1}^N,
          \quad  a_k  =  \left\{ \begin{array}{lll}
                                  0 \, & \mbox{ for }&
k> 0  \\
\displaystyle{\frac{\I}{ {\, 2 \,}}} \,
                                I_p
                                  & \mbox{ for }& k =
0  \\
                                \, \I \, I_p
                                  & \mbox{ for }& k <
0
                          \end{array} \right.
, \qquad J=\begin{bmatrix} 0 & I_p \\ I_p & 0 \end{bmatrix};
\\ &  \label{nu} 
\Phi_1(N) = \left[
\begin{array}{c}
I_{p}  \\ I_p \\ \cdots \\ I_{p}
\end{array}
\right], \quad \Phi_2(N) =  \left[
\begin{array}{l}
s_0/2  \\ s_0/2 + s_{-1} \\ \cdots \\ s_0/2+ s_{-1} + \ldots +
s_{1-N}
\end{array}
\right]
 +\I \Phi_1(N) \nu, 
\end{align}
and $\nu = \nu^*$. 
Moreover, in view of Verblunsky-type results \cite[Theorem 2.6]{ALS-JAT} (see also the proof of \cite[Theorem 2.4]{ALS-JAT})
there is a one to one correspondence between the sets $\{\nu\}\cup \{s_{-k}\}$ $(0 \leq k<\infty)$ such that \eqref{P5} holds
and Dirac systems \eqref{I1}, \eqref{I2}. This correspondence is given (in one direction) by \eqref{P4} and by the equalities
\begin{align}& \label{P7-}
\nu=\Im(\a_0), \quad s_0=\a_0+\a_0^*.
\end{align}
The transfer matrix function in Lev Sakhnovich form \cite{SaL1, SaLspth} is given by
the formula
\begin{align}& \label{P7}
w_A(N, \la)= I_{2p}-\I J\Pi(N)^*S(N)^{-1}\big(A(N)-\la I_{np}\big)^{-1}\Pi(N),
\end{align} 
and the recovery of Dirac system from the set $\{\nu\}\cup \{s_{-k}\}$ (the Verblunsky-type correspondence in another direction)
is equivalent to the factorisation of the transfer matrix functions $w_A(N, \la)$ \cite{ALS-JAT}. More precisely, by virtue of
\cite[(2.11)]{ALS-JAT} we have
\begin{align}& \label{P8}
W_{N}(\la)=\la^{-N}(\la +\I)^{N}K^*w_A(N,-\la/2)K \quad (0<N <\infty).
\end{align} 
\begin{Rk}\label{RkC+} If we substitute $\BC_+$ for $\BC_-$, Definition \ref{DnWD} defines Weyl functions in $\BC_+$.
According to \cite[Theorem 6.8]{FKRS08}, the Weyl functions $\vp(z)$ in $\BC_+$ and $\BC_-$ are unique
and are connected by the relation
\begin{align}& \label{P9}
\vp(\la)=\vp(\ov{\la})^*.
\end{align}
\end{Rk}
Relations \eqref{P4}, \eqref{P5} provide Verblunsky-type one to one correspondence between Dirac systems
and Toeplitz matrices (with additional matrix $\nu=\nu^*$ given by \eqref{P7-}) in one direction and  relations \eqref{P7}, \eqref{P8} provide Verblunsky-type mapping in the opposite
direction (as well as constitute the main part of the solution of the inverse problem to recover Dirac system from the Weyl function), see \cite{ALS-JAT}.

\section{Asymptotics of the analog \\ of reproducing kernel}\label{CD}
\setcounter{equation}{0}
{\bf 1.} Christoffel functions, reproducing kernels and Christoffel--Darboux  formula are important components
of the theory of orthogonal polynomials (see, e.g., \cite{BreLaS, Lub, Ne, NeTo} and references therein). In this section, we consider an analog of the Christoffel--Darboux  formula
for the discrete Dirac system and study the  asymptotics of the analog of reproducing kernel. First, we prove the following Christoffel--Darboux-type  formula.
\begin{Tm}\label{TmCD} Let $W_k(\la)$ be the fundamental solution of the discrete Dirac system \eqref{I1}, \eqref{I2}
normalised by \eqref{P1}. Then, we have
\begin{align} &  \label{C1}
\sum_{k=0}^{N}c(\la,\mu)^kW_k(\ov{\mu})^*C_kW_k(\la)=\I \frac{1+\la \mu}{\mu - \la}\Big(c(\la,\mu)^{N+1}W_{N+1}(\ov{\mu})^*jW_{N+1}(\la)-j\Big),
\end{align} 
where
\begin{align}& \label{C2}
c(\la,\mu)=\frac{\la \mu}{1+\la \mu}.
\end{align} 
\end{Tm}
We note that $c(\la, \ov{\la})=q(\la)$ for $q$ introduced in \eqref{P3}.
\begin{proof} of Theorem \ref{TmCD}. Taking into account \eqref{I1} and \eqref{I2}, we obtain
\begin{align}\nn
W_{k+1}(\ov{\mu})^*jW_{k+1}(\la)&=W_{k+1}(\ov{\mu})^*\Big(I_{2p}+\frac{\I}{\mu}C_k j\Big)j\Big(I_{2p}-\frac{\I}{\la}jC_k\Big)
W_{k+1}(\la)
\\ &  \label{C3}
=W_{k}(\ov{\mu})^*\Big(\frac{1}{c(\la, \mu)}j+\frac{\I(\la - \mu)}{\la \mu} C_k    \Big)W_{k}(\la).
\end{align} 
Formula \eqref{C3} yields
\begin{align}& \label{C4}
W_{k}(\ov{\mu})^* C_k    W_{k}(\la)=\frac{\I(1+\la \mu)}{\mu-\la}\Big(c(\la,\mu)W_{k+1}(\ov{\mu})^*jW_{k+1}(\la)-W_{k}(\ov{\mu})^*jW_{k}(\la)\Big).
\end{align} 
Substitute \eqref{C4} into the left-hand side of \eqref{C1} in order to derive \eqref{C1}.
\end{proof}
{\bf 2.} Let us fix the set $\{\nu\}\cup \{s_{-k}\}$ $(0 \leq k<\infty)$ such that \eqref{P5} holds.
For asymptotical results in this section, we use \cite[Example 1]{ALS87} and \cite[Theorem 6.5]{FKRS08}.
The notations in \cite{ALS87} differ from the notations here and we give some explanations.
We recall that the Weyl function $\vp(\la)$ is constructed in \cite[Theorem 6.5]{FKRS08} as the unique function
belonging to  the intersection of the Weyl circles generated by some linear fractional transformations.
The corresponding linear fractional transformations are given by  \cite[(5.7)]{FKRS08}
and closely related transformations are given by the formulas \cite[(11)]{ALS87}.
The matrices of coefficients of these transformations  are denoted by $\clw$ in \cite{FKRS08} and by $\mathfrak{A}$ in \cite{ALS87},
where  
\begin{align}& \label{C5!}
\clw(\la)=\clw_N(\la)=KW_N(\ov{\la})^* \quad (N\in \BN),
\end{align} 
$\mathfrak{A}$ (in the notations of this paper) takes the form
\begin{align}& \label{C11}
\mathfrak{A}(\zeta)=\mathfrak{A}_N(\zeta)=j\big(I_{2p}+\I\zeta J \Pi(N)^*\big(I_{pN}+\zeta A(N)^*\big)^{-1}S(N)^{-1}\Pi(N)\big)j,
\end{align} 
and $K$ and $J$ are introduced in \eqref{P3} and \eqref{P6+}, respectively. 
Simple calculations show that these matrices (matrix functions) satisfy the equality
\begin{align}& \label{C5-}
\clw(\la)=\left(\frac{\la -\I}{\la}\right)^N j J{\mathfrak{A}}\left(\frac{2}{\la}\right)JjK.
\end{align} 
Moreover, it is easy to see that 
\begin{align} &  \label{C13}
 K^*JK=j, \quad jJK=K^*.
\end{align} 

We note that linear-fractional transformations considered in \cite{ALS87} have the form
\begin{align} &  \label{LF1}
 \om_N(\zeta)=\I\big(a(\zeta)R(\zeta)+b(\zeta)Q(\zeta)\big)\big(c(\zeta)R(\zeta)+d(\zeta)Q(\zeta)\big)^{-1}, \quad \mf{A}_N=:\begin{bmatrix} a & b \\c & d\end{bmatrix},
 \end{align} 
 where $\zeta \in \BC_+$; $a,b,c,d$ are $p\times p $ blocks of $\mf{A}$ and $\{R(\zeta), \, Q(\zeta)\}$ are pairs of  $p \times p$ matrix-valued functions,
 which are meromorphic in $\BC_+$ and satisfy inequalities:
 \begin{align} &  \label{LF2}
R(\zeta)^*R(\zeta)+Q(\zeta)^*Q(\zeta)>0, \quad \begin{bmatrix}R(\zeta)^* & Q(\zeta)^* \end{bmatrix}J \begin{bmatrix}R(\zeta) \\ Q(\zeta) \end{bmatrix}\geq 0
 \end{align} 
(excluding, possibly, some isolated points). Such pairs $\{R(\zeta), \, Q(\zeta)\}$  are called nonsingular and $J$-nonnegative.
The first equality in \eqref{C13} shows that the matrix $JjK$ maps $-j$-nonegative pairs $\{\wt R,\wt Q\}$ (more precisely, the
columns col$\begin{bmatrix}\wt R & \wt Q\end{bmatrix}$) which are used in the linear fractional
transformations \cite[(5.7)]{FKRS08} onto the set of $J$-nonegative pairs $\{R,Q\}$ which are used in the linear fractional
transformations  \cite[(11)]{ALS87}.
In view of \cite[(5.7)]{FKRS08}, \cite[(11)]{ALS87} and the equality \eqref{C5-} above,
it follows that our Weyl function $\vp$ is connected with the function $\om$ corresponding to $\{\nu\}\cup \{s_{-k}\}$ $(0 \leq k<\infty)$
in \cite{ALS87}  by the relation
\begin{align}& \label{C5}
\om(\zeta )=-\vp(2/\zeta ).
\end{align} 
Since $\vp(\la)$ is Herglotz function  in $\BC_-$, $\om(\zeta )$ belongs to Herglotz class in $\BC_+$ (i.e., $\Im\big(\om(\zeta )\big)\geq 0$).

{\bf 3.}
Now, consider Herglotz representation
\begin{align}& \label{C6}
\om(\zeta )=\b \zeta +\g+\int_{-\infty}^{\infty}\frac{1+t\zeta }{(t-\zeta )(1+t^2)}d\tau(t); \\
& \label{C6'}
 \b\geq 0, \quad \g=\g^*, \quad \int_{-\infty}^{\infty}{(1+t^2)^{-1}}{d\tau(t)}<\infty,
\end{align} 
and assume that  $\tau^{\prime}$ (the positive semi-definite derivative of the absolutely continuous part of $\tau$) satisfies the Szeg{\H o} condition
\begin{align}& \label{C7}
\int_{-\infty}^{\infty}(1+t^2)^{-1}\ln\big(\det\tau^{\prime}(t)\big)>-\infty.
\end{align} 
The factors of $\tau^{\prime}$  will play a crucial role in our further considerations. The factorisation  of positive semi-definite integrable matrix functions
is one of the classical domains connected with the names of  A. Beurling, N. Wiener, P.R.~Masani, H. Helson, D. Lowdenslager, M.G. Krein, Yu.A. Rozanov,
D.~Sarason and many others (see, e.g., \cite{HeLo, MaWi, RoRo, Roz, Wien} and the bibliography in the interesting paper \cite{KiKa}).
We will use the results formulated in \cite[Theorem 9]{HeLo} and \cite[Theorem 5.2]{KiKa}.

Notice that the functions $(\overline{\zeta_0}z -\zeta_0 )(z -1)^{-1}$,
where $\zeta_0\in \BC_+$, map unit disk $\BD$ ($|z|<1$) onto $\BC_+$. Let us introduce the {\it class of functions $\wt H$ on $\BC_+$}.
We say that the $p \times p$ matrix function $G(\zeta)$ belongs to $\wt H$ if 
the entries of $G\left(\frac{\overline{\zeta_0}z -\zeta_0 }{z -1}\right)$ 
belong to the Hardy class $H^2(\BD)$ and $\det \left(G\left(\frac{\overline{\zeta_0}z -\zeta_0 }{z -1}\right)\right)$ is an outer function.
Setting
\begin{align}& \nn
t=(\overline{\zeta_0}z -\zeta_0 )(z -1)^{-1}=\overline{\zeta_0}+(\overline{\zeta_0}-{\zeta_0})(z-1)^{-1}, \quad z=\E^{\I\theta} \quad( 0\leq \theta <2\pi),
\end{align} 
and taking into account that $\ov{z}(z-1)=-\ov{(z-1)}$ and $\ov{z}(\overline{\zeta_0}z -\zeta_0 )=-\ov{\big(\overline{\zeta_0}z -\zeta_0 \big)}$ we obtain
\begin{align}& \label{C7+}
\frac{dt}{1+t^2}=-\frac{\I(\overline{\zeta_0}-{\zeta_0})zd\t}{(z-1)^2+(\overline{\zeta_0}z -\zeta_0 )^2}=
\frac{\I(\overline{\zeta_0}-{\zeta_0})d\t}{|z-1|^2+|\overline{\zeta_0}z -\zeta_0 |^2}.
\end{align} 
From \eqref{C6'} and \eqref{C7+}, it follows that 
$$\int_0^{2\pi}\tau^{\prime}\left(\frac{\overline{\zeta_0}\E^{\I \t} -\zeta_0 }{\E^{\I \t} -1}\right)d\t<\infty.$$
In other words, $\tau^{\prime}\left(\frac{\overline{\zeta_0}\E^{\I \t} -\zeta_0 }{\E^{\I \t} -1}\right)$ is integrable.
From  \eqref{C7} and \eqref{C7+}  we derive that $\int_0^{2\pi}\ln\det\left(\tau^{\prime}\left(\frac{\overline{\zeta_0}\E^{\I \t} -\zeta_0 }{\E^{\I \t} -1}\right)\right)d\t>-\infty$,
and (in view of the equality $\ln(\det  (\tau^{\prime}))=p\, {\rm tr}(\ln (\tau^{\prime}))$) the condition (74) of   \cite[Theorem 9]{HeLo} is fulfilled.
Thus, there is a factorisation
\begin{align}& \label{C8}
\tau^{\prime}(t)=G_{\tau}(t)^*G_{\tau}(t),
\end{align} 
where $G_{\tau}(\zeta )\in \wt H$ and $G_{\tau}(t)$ is the boundary function of  $G_{\tau}(\zeta )$. Here, we reversed the order of factors in
\cite[Theorem 9]{HeLo}, which does not matter (see, e.g., \cite[p. 195]{HeLo}).  One can use \cite[Theorem 5.2]{KiKa} on inner-outer factorisation
to show that    an outer matrix function $G\left(\frac{\overline{\zeta_0}z -\zeta_0 }{z -1}\right)$ (such that $\det(G)$ is a scalar outer function)
may be chosen uniquely up to a constant unitary factor. 
The definition of $\wt H$ does not depend on the choice of $\zeta_0\in \BC_+$ because the Hardy classes $H^p$ are invariant under conformal one-to-one transformations
of $\BD$.

 The matrix functions $\clr_k(\la,\mu)\,$ (i.e., $\rho_k(\la,\mu)$ in \cite{ALS87})
are introduced by the equality
\begin{align}& \label{C9}
\clr_k(\la,\mu)=\Phi_1(k)^*\big(I_{kp}+\la A(k)^*\big)^{-1}S(k)^{-1}\big(I_{kp}+\mu A(k)\big)^{-1}\Phi_1(k).
\end{align}
\begin{Rk}\label{RkNond}
Clearly,  the matrices $\clr_k(\zeta,\ov{\zeta})$ are well-defined
and positive-definite for $\zeta \not= -2 \I$. According to \cite[Theorem 2]{ALS87}, the sequence $\clr_k(\zeta,\ov{\zeta})$ $(\zeta\in \BC, \,\,\zeta \not=- 2\I)$ is nondecreasing.
\end{Rk} 
From \cite[Example 1]{ALS87}, \cite[Theorem 4]{ALS87} and \cite[Remark 1]{ALS87}, we obtain the following theorem.
\begin{Tm}\label{TmRho} Let us fix the set of $p\times p$ matrices $\{\nu\}\cup \{s_{-k}\}$ $(0 \leq k<\infty)$, where $\nu=\nu^*$, $s_0=s_0^*$ and we set $s_k=s_{-k}^*$
for $k>0$. Assume that the inequalities  \eqref{P5} hold for all $1 \leq N <\infty$ and that Szeg{\H o} condition \eqref{C7} is fulfilled.
Then, 
\begin{align}& \label{C10}
\lim_{k\to \infty}\clr_k(\zeta,\ov{\zeta})^{-1}=2\pi\I (\ov{\zeta}-\zeta)G_{\tau}(\zeta)^*G_{\tau}(\zeta) \qquad (\zeta \in \BC_+).
\end{align} 
\end{Tm}
\begin{Rk}\label{RkUn} In view of Remark \ref{RkNond}, the limit in \eqref{C10} is uniform for $\zeta$ belonging to the compact subsets
of $\BC_+$.
\end{Rk}
{\bf 4.} {\it The analogs of the reproducing kernels} are matrix functions $W_k(\ov{\mu})^*jW_k(\la)$ (see, e.g., \eqref{C4}).
Let us rewrite this expression in terms of $\mf{A}_k$. From the definitions \eqref{P3}, \eqref{P7} and \eqref{C11} we derive
\begin{align}& \label{C12}
w_A(k,-\ov{\mu}/2)^*=jJ\mf{A}_k(2/\mu)Jj, \quad w_A(k,-{\la}/2)=jJ\mf{A}_k(2/\ov{\la})^*Jj.
\end{align} 
Using \eqref{P8}, \eqref{C13} and \eqref{C12} we obtain
\begin{align}\nn
W_k(\ov{\mu})^*jW_k(\la)&=\frac{(\mu-\I)^k(\la+\I)^k}{(\mu \la)^k} K^*jJ\mf{A}_k(2/\mu)JjKjK^*jJ\mf{A}_k(2/\ov{\la})^*JjK
\\ & \label{C14}
= -\frac{(\mu-\I)^k(\la+\I)^k}{(\mu \la)^k}K\mf{A}_k(2/\mu)J \mf{A}_k(2/\ov{\la})^*K^*.
\end{align} 
Thus, we can study the asymptotics of $\mf{A}_k(\zeta)J \mf{A}_k(\ov{\xi})^*$ instead of the asymptotics of  $W_k(\ov{\la})^*jW_k(\mu)$.
We set
\begin{align}& \label{C14'}
\clm(k,\zeta,\xi):=\mf{A}_k(\zeta)J \mf{A}_k(\ov{\xi})^*.
\end{align} 

Relations \eqref{C12} and \eqref{C14'} yield
\begin{align}& \label{C15}
\clm(k,\zeta,\xi)=-Jjw_A(k,-1/\ov{\zeta})^*Jw_A(k,-1/{\xi})jJ.
\end{align} 
Hence, using either the properties of the transfer function $w_A$ (see, e.g., \cite[Corollary 1.15]{SaSaR}) or  direct calculation (which takes into account \eqref{P6}),
we have
\begin{align} \nn 
\clm(k,\zeta,\xi)=&J+\I(\xi-\zeta)Jj
\Pi(k)^*
\\ & \label{C16}
\times (I_{kp}+\zeta A(k)^*)^{-1}S(k)^{-1}(I_{kp}+\xi A(k))^{-1}\Pi(k)
jJ.
\end{align} 
Partition $\clm$ into $p\times p$ blocks $\clm=\{\clm_{ik}\}_{i,k=1}^2$. From \eqref{C9}
and \eqref{C16} it follows that
\begin{align}& \label{C17}
\clm_{22}(k,\zeta,\xi)=\I (\xi-\zeta)\clr_k(\zeta, \xi).
\end{align} 

Thus, Theorem \ref{TmRho} describes the asymptotics of $\clm_{22}(k,\zeta,\ov{\zeta})$:
\begin{align}& \label{C10'}
\lim_{k\to \infty}\clm_{22}(k,\zeta,\ov{\zeta})=(1/2\pi) \big(G_{\tau}(\zeta)^*G_{\tau}(\zeta)\big)^{-1} \qquad (\zeta \in \BC_+).
\end{align}
Another way to obtain this asymptotics is to use the note \cite{AKr}. According to Remark \ref{RkUn}, the limit \eqref{C10'} is uniform.

{\bf 5.} Next, using (for the asymptotics of $\clm_{22}(k,\zeta,\ov{\xi})$) 
an approach from the theory of orthogonal
polynomials (see, e.g., \cite{ApNi}) we prove
the following theorem.
\begin{Tm} \label{TmM} Let us fix the set of $p\times p$ matrices $\{\nu\}\cup \{s_{-k}\}$ $(0 \leq k<\infty)$, where $\nu=\nu^*$, $s_0=s_0^*$ and we set $s_k=s_{-k}^*$
for $k>0$. Assume that the inequalities  \eqref{P5} hold for all $1 \leq N <\infty$ and that Szeg{\H o} condition \eqref{C7} is fulfilled. Let  complex values
$\zeta$ and $\xi$ belong to some compact subset in $\BC_+\backslash\{2\I\}$.
Then, uniformly with respect to $\zeta$ and $\xi$ in this compact subset, we have
\begin{align}& \label{C30}
\lim_{k\to \infty}\clm(k,\zeta, \ov{\xi})=\frac{1}{2\pi}\, \left[\begin{array}{c}-\I\om ({\zeta} ) \\ I_p
\end{array}\right]\, {G}_{\tau}({\zeta} )^{-1}\, \left({G}_{\tau}(\xi )^*
\right)^{-1}\,\left[ \begin{array}{lr} \I \om(\xi )^* &
I_p\end{array}\right] ,
\end{align} 
where $\clm(k,\zeta,\ov{\xi})=\mf{A}_k(\zeta)J \mf{A}_k({\xi})^*$, $\mf{A}_k$ is given by \eqref{C11}, $\om$ is given by \eqref{C5}, 
and $G$ is the factoring multiplier from \eqref{C8} $(G(\ze)\in \wt H)$.
\end{Tm}
\begin{proof}. Step 1.
Denote by $\Gamma_r$ the curve $\big|\left({\zeta}
-\xi \right)\left({\zeta} -\overline{\xi}\right)^{-1}\big| =r\leq 1$, where $\xi$ is some fixed point in $\BC_+$,
choose anticlockwise orientation for this $\Gamma_r$ and put
\begin{align} \nn &
\psi (k,r, \xi):=\frac{1}{2\pi \I}\int_{\Gamma_r}\left(
2\pi{G}_{\tau} (\xi )\clm_{22}(k,{\zeta} ,\ov{\xi} )^*{G}_{\tau}({\zeta} )^*-I_p\right)
\\ \label{C20} &
{\hspace{2cm}} \times \left( 2\pi{G}_{\tau} ({\zeta} )\clm_{22}(k,{\zeta} , \ov{\xi} )
{G}_{\tau}(\xi )^*-I_p\right)\, \frac{\left(\xi -\overline{\xi}
\right) \, d{\zeta}}{\left({\zeta} -\xi\right)\left({\zeta}
-\overline{\xi}\right)} ;
\\ \nn &
\widetilde{\psi}(k,r,\xi):= 2\pi \I \int_{\Gamma_{r}}{G}_{\tau} (\xi
) \clm_{22}(k,{\zeta} , \ov{\xi} )^* {G}_{\tau}({\zeta} )^*{G}_{\tau} ({\zeta})\clm_{22} (k,{\zeta} ,\ov{\xi} )
\\ \label{C21} &
{\hspace{3cm}}
\times \frac{{G}_{\tau}(\xi )^*\left(\overline{\xi}-\xi \right) 
d{\zeta}}{\left({\zeta} -\xi\right)\left({\zeta} -\overline{\xi}\right)}  .
\end{align}
Clearly, for $r<1$ we have
\begin{align}& \label{C23}
\frac{1}{2\pi \I}\int_{\Gamma_r}{G}_{\tau} ({\zeta} )\clm_{22}(k,{\zeta} , \ov{\xi} )
{G}_{\tau}(\xi )^* \frac{\left(\xi -\overline{\xi}
\right)\, d{\zeta}}{\left({\zeta} -\xi\right)\left({\zeta}
-\overline{\xi}\right)}={G}_{\tau} (\xi )\clm_{22}
(k,\xi ,\ov{\xi})\, {G}_{\tau}(\xi )^*.
\end{align}
Moreover, the curve $\Gamma_r$ may be rewritten in the form
$\left({\zeta}
-\xi \right)\left({\zeta} -\overline{\xi}\right)^{-1}=r\E^{\I\t}$,
which yields:
$$\zeta-\ov{\xi}=\frac{\xi -\ov{\xi}}{1-r \E^{\I\t}}, \qquad \zeta - \xi=(\xi - \ov{\xi})\frac{r\E^{\I\t}}{1-r\E^{\I\t}}.$$
Hence, the equality
\begin{align}& \label{C23+}
\frac{\left(\xi -\overline{\xi}
\right)\, d{\zeta}}{\left({\zeta} -\xi\right)\left({\zeta}
-\overline{\xi}\right)}=\I d \t
\end{align}
follows, and (in view of \eqref{C23}) we obtain
\begin{align}&\nn
\frac{1}{2\pi \I}\int_{\Gamma_r}
{G}_{\tau} (\xi )\clm_{22}(k,{\zeta} ,\ov{\xi} )^*{G}_{\tau}({\zeta} )^*
\frac{\left(\xi -\overline{\xi}
\right)\, d{\zeta}}{\left({\zeta} -\xi\right)\left({\zeta}
-\overline{\xi}\right)}
\\ \nn &
=\frac{1}{2\pi}\int_0^{2\pi}
{G}_{\tau} (\xi )\clm_{22}(k,{\zeta(\t)} ,\ov{\xi} )^*{G}_{\tau}({\zeta(\t)} )^*d\t
\\  &  \nn
=\left(\frac{1}{2\pi}\int_0^{2\pi}{G}_{\tau} ({\zeta(\t)} )\clm_{22}(k,{\zeta(\t)} , \ov{\xi} )
{G}_{\tau}(\xi )^*d\t\right)^*
\\  & \label{C24}
=\big({G}_{\tau} (\xi )\clm_{22}
(k,\xi ,\ov{\xi})\, {G}_{\tau}(\xi )^*\big)^*
={G}_{\tau} (\xi )\clm_{22}
(k,\xi ,\ov{\xi})\, {G}_{\tau}(\xi )^*.
\end{align}
Here, we used the equality $\clm_{22}
(k,\xi ,\ov{\xi})=\clm_{22}
(k,\xi ,\ov{\xi})^*$, which is immediate from \eqref{C15}. Equalities \eqref{C20}--\eqref{C23} and \eqref{C24} imply that
\begin{equation}\label{C22}
\psi (k,r,\xi)=I_p+\widetilde{\psi}(k,r,\xi)- 4\pi{G}_{\tau} (\xi )\clm_{22}
(k,\xi ,\ov{\xi})\, {G}_{\tau}(\xi )^*.
\end{equation}

We note that $\Gamma_1=\BR$ (where the curves $\Gamma_r$ were introduced at the beginning of the proof),
and so (in view of  \eqref{C8}, \eqref{C9}, \eqref{C17} and \eqref{C21}) we have
\begin{align}\nn 
 \widetilde{\psi}(k,1,\xi)=&2\pi \I (\ov{\xi}-\xi)G_{\tau}(\xi)\Phi_1(k)^*(I_{kp}+{\xi}A(k)^*)^{-1}S(k)^{-1}
 \\ \nn &
\times \int_{-\infty}^{\infty}(I_{kp}+{\zeta} A(k))^{-1}\Phi_1(k)\tau^{\prime}(\zeta)\Phi_1(k)^*(I_{kp}+\zeta A(k)^*)^{-1}d\zeta
\\ \label{C26} &
\times S(k)^{-1}(I_{kp}+\ov{\xi} A(k))^{-1}\Phi_1(k)G_{\tau}(\xi)^*.
\end{align}
From \cite[(7) and (9)]{ALS87}, we obtain the following representation of $S(k)$ (in the present notations):
\begin{align}\nn 
S(k)=&A^{-1} \Phi_1(k)\b \Phi_1(k)^* (A^{-1})^*   
\\ \label{C25} &
  +\int_{-\infty}^{\infty}(I_{kp}+t A(k))^{-1}\Phi_1(k)d\tau(t)\Phi_1(k)^*(I_{kp}+t A(k)^*)^{-1},
\end{align}
where $\Phi_1$ is introduced in \eqref{nu} and $\b$ and $\tau$ are  given by the Herglotz representation \eqref{C6}.
Thus, it is easy to see that
\begin{align}\label{C25+} &
 \int_{-\infty}^{\infty}(I_{kp}+{\zeta} A(k))^{-1}\Phi_1(k)\tau^{\prime}(\zeta)\Phi_1(k)^*(I_{kp}+\zeta A(k)^*)^{-1}d\zeta
\leq S(k).
\end{align}
Finally, relations \eqref{C17}, \eqref{C26} and \eqref{C25+} yield
\begin{equation}\label{C27}
\widetilde{\psi }(k,1,\xi) \leq 2\pi\,{G}_{\tau} (\xi )\,\clm_{22} (k,{\xi}
, \ov{\xi} ) \, {G}_{\tau}(\xi )^*.
\end{equation}
  This implies that the entries of ${G}_{\tau} ({\zeta} )\clm_{22} ({\zeta} , \ov{\xi} )$  belong (after substitution
  ${\zeta} =
\left(\, \overline{\zeta_0}\, z -\zeta_0\right)\left( z -1 \right)^{-1},$ $\zeta_0\in \BC_+$)
to the space
$H^2$ of the analytic functions of $z$ in the unit disk. Then, using theorem of F.~Riesz  (see, e.g., formula
(4.1.1) in \cite{Pri}) and properties of subharmonic functions we
have
\begin{equation}\label{C28}
\widetilde{\psi}(k,r,\xi)\leq \widetilde{\psi}(k,1,\xi) \quad (r\leq 1).
\end{equation}
According to \eqref{C20} and \eqref{C23+} the inequality $\psi (k,r,\xi)\geq 0$ is valid.
Since $\psi (k,r,\xi)\geq 0$, we take into account \eqref{C10'}, \eqref{C27} and \eqref{C28} and derive  from
\eqref{C22}  that
\begin{equation}\label{C29}
\lim\limits_{k\to\infty}\psi (k,r,\xi)=0,
\end{equation}
uniformly in $\xi$ and $r$. Using \eqref{C29}
and expanding 
$${G}_{\tau} \left( \left(\, \overline{\zeta_0}\,z -\zeta_0\right)\left( z -1 \right)^{-1}\right)
\clm_{22} \left(\left(\,\overline{\zeta_0}\, z -\zeta_0\right)\left( z -1 \right)^{-1} , \ov{\xi} \right)$$
 in series in $z $,  we obtain the assertion:
 \begin{align}& \label{C10+}
\lim_{k\to \infty}\clm_{22}(k,\zeta,\ov{\xi})=(1/2\pi) \big(G_{\tau}(\xi)^*G_{\tau}(\zeta)\big)^{-1} \qquad (\zeta,\xi \in \BC_+)
\end{align}
uniformly on the compacts in $\BC_+$. Formula \eqref{C10+} coincides with the restriction of \eqref{C30} to $\clm_{22}(k,\zeta, \ov{\xi})$.

Step 2. The set of the matrix-valued functions given by the linear-fractional (M\"obius) transformations \eqref{LF1}, where $\{R(\zeta), \, Q(\zeta)\}$  are nonsingular, $J$-nonnegative
pairs, is denoted by $\cln(\mf{A}_N)$ $(N\in \BN)$, and the set of values which these matrix functions take at $\zeta \in \BC_+$ is denoted by $\cln(\mf{A}_N)(\zeta)$. 
The sets $\cln(\mf{A}_N)$ are embedded (i.e., $\cln(\mf{A}_N)\subseteq \cln(\mf{A}_{\wh N})$ for $N>\wh N$), and their intersection consists of one function
$\om(\zeta)$; see, for instance, Proposition 5.7 and Theorem 6.4 in \cite{FKRS08}.  (One easily deletes the requirement from \cite[Definition 5.3]{FKRS08}
that the pairs $\{R,Q\}$ are well-defined at some fixed point). Moreover, according to \eqref{C5} and \cite[p. 227]{FKRS08} we have
\begin{align}\label{C41}
\bigcap_{N\geq 1}\cln(\mf{A}_N)(\zeta)=\{\om(\zeta)\}.
\end{align}
Next, we show that $\cln(\mf{A}_N)(\zeta)$ form Weyl-type disks and consider these disks.
By virtue of \eqref{C14'} and \eqref{C16} (see also \cite[(17)]{ALS87}), we have
\begin{align}& \label{C31}
\mf{A}_N(\zeta)J\mf{A}_N(\ov{\zeta})^*=J=\mf{A}_N(\ov{\zeta})^*J\mf{A}_N(\zeta).
\end{align}
Since $\mf{A}_N(\ov{\zeta})^*J\mf{A}_N(\zeta)=J$, it is immediate that equality \eqref{LF1} (where \eqref{LF2} holds) is equivalent to
\begin{align}& \label{C32}
\begin{bmatrix}\I \om_N(\zeta)^* & I_p \end{bmatrix}
J\mf{A}_N(\ov{\zeta})J\mf{A}_N(\ov{\zeta})^*J \begin{bmatrix}-\I \om_N(\zeta) \\ I_p \end{bmatrix} \geq 0.
\end{align}
Setting 
\begin{align}& \label{C33}
\mf{F}(\zeta):=J\mf{A}_N(\ov{\zeta})J\mf{A}_N(\ov{\zeta})^*J=\{\mf{F}_{ik}(\zeta)\}_{i,k=1}^2,
\end{align}
where $\mf{F}_{ik}$ are $p\times p$ blocks of $\mf{F}$, we rewrite \eqref{C32} in the form
$$\om_N(\zeta)^*\big(-\mf{F}_{11}(\zeta)\big)\om_N(\zeta)+\I\big(\mf{F}_{21}(\zeta)\om_N(\zeta)-\om_N(\zeta)^*\mf{F}_{12}(\zeta)\big)\leq \mf{F}_{22}(\zeta).$$
Equivalently, we have
\begin{align}& \nn
\om_N(\zeta)^*\big(-\mf{F}_{11}(\zeta)\big)\om_N(\zeta)+\I\big(\mf{F}_{21}(\zeta)\om_N(\zeta)-\om_N(\zeta)^*\mf{F}_{12}(\zeta)\big)
\\ & \label{C34}
+
\mf{F}_{21}(\zeta)\big(-\mf{F}_{11}(\zeta)\big)^{-1}\mf{F}_{12}(\zeta)
\leq \mf{F}_{22}(\zeta)-
\mf{F}_{21}(\zeta)\mf{F}_{11}(\zeta)^{-1}\mf{F}_{12}(\zeta),
\end{align}
where 
\begin{align}& \label{C34+}
-\mf{F}_{11}(\zeta)=\I(\ov{\zeta}-\zeta)\clr_N(\ov{\zeta},\zeta)>0.
\end{align}
From \eqref{C33} (using again \eqref{C31}), we obtain
\begin{align}& \label{C35}
\mf{F}(\zeta)^{-1}=\big\{\big(\mf{F}(\zeta)^{-1}\big)_{ik}\big\}_{i,k=1}^2=\mf{A}_N(\zeta)J\mf{A}_N(\zeta)^*,
\end{align}
where $\big(\mf{F}(\zeta)^{-1}\big)_{ik}$ are $p\times p$ blocks of $\mf{F}(\zeta)^{-1}$.
Taking into account \eqref{C14'}, \eqref{C17} and \eqref{C35} we derive
\begin{align}& \label{C36}
\big(\mf{F}(\zeta)^{-1}\big)_{22}=\I (\ov{\zeta}-\zeta)\clr_N(\zeta, \ov{\zeta})>0.
\end{align}
Since $\big(\mf{F}(\zeta)^{-1}\big)_{22}$ is invertible, we express $\big(\big(\mf{F}(\zeta)^{-1}\big)_{22})^{-1}$
in terms of the blocks $\mf{F}_{ik}$ (see, e.g., \cite[Section 1.2]{SaLspth}):
\begin{align}& \label{C37}
\big(\big(\mf{F}(\zeta)^{-1}\big)_{22})^{-1}=\mf{F}_{22}(\zeta)-
\mf{F}_{21}(\zeta)\mf{F}_{11}(\zeta)^{-1}\mf{F}_{12}(\zeta).
\end{align}
In view of \eqref{C34+} and \eqref{C36}, the  matrix functions $\Lam_l(N,\zeta)$ and $\Lam_r(N,\zeta)$
(which we will show to be the left and right radii)
are well-defined by the relations
\begin{align}& \label{C38}
\Lam_l(N,\zeta)^2=\big(-\mf{F}_{11}(\zeta)\big)^{-1}, \quad \Lam_l(N,\zeta)>0; 
\\ & \label{C39}
\Lam_r(N,\zeta)^2=\big(\big(\mf{F}(\zeta)^{-1}\big)_{22})^{-1}, \quad \Lam_r(N,\zeta)>0.
\end{align}
Using \eqref{C37}--\eqref{C39}, we rewrite \eqref{C34} in the form 
$$(\Lam_l^{-1}\om_N-\I\Lam_l\mf{F}_{12})^*(\Lam_l^{-1}\om_N-\I\Lam_l\mf{F}_{12})\leq 
\Lam_r^2.$$
That is, we  parametrize $\cln(\mf{A}_N)(\zeta)$ via contractive $p\times p$ matrices $u$ as a matrix circle
\begin{align}& \label{C40}
\om_N(\zeta)=\Lam_l(N,\zeta)u\Lam_r(N,\zeta)+\I \Lam_l(N,\zeta)^2\mf{F}_{12}(\zeta) \qquad (u^*u\leq I_p),
\end{align}
where $\Lam_l(N,\zeta)$ and $\Lam_r(N,\zeta)$ are, indeed, the left and right radii,  and 
$N$ is omitted in the notations connected with $\mf{F}$. According to \eqref{C34+}, \eqref{C38} and \eqref{C36}, \eqref{C39},
we obtain
\begin{align}& \label{C45}
\Lam_l(N,\zeta)^2=\I(\zeta -\ov{\zeta})^{-1}R_N\big(\ov{\zeta}, {\zeta}\big)^{-1}, 
\quad \Lam_r(N,\zeta)^2=\I(\zeta -\ov{\zeta})^{-1}R_N\big(\zeta,\ov{\zeta}\big)^{-1}.
\end{align}
Relations \eqref{C10}, \eqref{C45} and Remark \ref{RkUn}
imply that {\it uniformly on each compact in $\BC$} we have
\begin{align}& \label{C42}
\lim_{k\to \infty}\Lam_r(k,\zeta)=\sqrt{ 2\pi G_{\tau}(\zeta)^*G_{\tau}(\zeta)}>0.
\end{align}
Hence, equalities \eqref{C41} and \eqref{C40} yield
\begin{align}& \label{C43}
\lim_{k\to \infty}\Lam_l(k,\zeta)=0 \quad  (\zeta \in \BC_+, \,\, \zeta\not=2\I).
\end{align}
Taking into account \eqref{C34+}, \eqref{C38} and \eqref{C43} we obtain
\begin{align}& \label{C44}
\lim_{k\to \infty}\clr_k(\ov{\zeta},\zeta)^{-1}=0 \quad (\zeta \in \BC_+, \,\, \zeta\not=2\I).
\end{align}
Moreover, according to \cite[Theorem 2]{ALS87} the sequences of matrices $\clr_k(\ov{\zeta}, \zeta)$ (matrices $\rho_k$ in the notations of 
\cite{ALS87}) are nondecreasing. Thus, \eqref{C44} (and so \eqref{C42} as well) holds {\it uniformly for the values $\zeta$ on any compact} $\clc \subset \BC_+ \backslash \{2\I\}$.

Step 3.   Using \eqref{C14'}, \eqref{C33} and \eqref{C35} (and substituting $\ov{\zeta}$ instead of $\zeta$), we rewrite \eqref{C37} in the form
\begin{align}& \label{C47}
\clm_{22}\big(N, \ov{\zeta}, {\zeta}\big)^{-1}=\mf{F}_{22}(\ov{\zeta})-
\mf{F}_{21}(\ov{\zeta})\clm_{22}\big(N,  {\zeta}, \ov{\zeta}\big)^{-1}\mf{F}_{12}(\ov{\zeta}) \quad (\zeta \not=\ov{\zeta}, \,\, \zeta \not=\pm 2\I).
\end{align}
In order to right down the necessary expressions for $\mf{F}_{22}, \, \mf{F}_{21}$ and $\mf{F}_{12}$, we need the representation 
\begin{align}\nn &
\begin{bmatrix} a(N,\zeta) & b(N,\zeta)\end{bmatrix}=a(N,\zeta)c(N,\zeta)^{-1}\begin{bmatrix} c(N,\zeta) & d(N,\zeta)\end{bmatrix}
\\  &  \label{C48} \hspace{9em}
-q(N,\zeta)\begin{bmatrix} 0 & d(N,\zeta)\end{bmatrix}, \quad 
\\  &  \label{C49}
q(N,\zeta):=a(N,\zeta)c(N,\zeta)^{-1}-b(N,\zeta)d(N,\zeta)^{-1}.
\end{align}
From \eqref{C14'}, \eqref{C33}, \eqref{C47} and \eqref{C48} (omitting variables $\zeta$ and $N$ in the matrix functions
$a,b,c,d$ and $q$), we obtain
\begin{align}\nn
\clm_{22}\big(N, \ov{\zeta}, {\zeta}\big)^{-1}=&ab^*+ba^*-\big(ac^{-1}\clm_{22}\big(N,  {\zeta}, \ov{\zeta}\big)-qdc^*\big)\clm_{22}\big(N,  {\zeta}, \ov{\zeta}\big)^{-1}
\\ & \label{C50} \times
\big(ac^{-1}\clm_{22}\big(N,  {\zeta}, \ov{\zeta}\big)-qdc^*\big)^*,
\end{align}
where $\clm_{22}\big(N,  {\zeta}, \ov{\zeta}\big)=cd^*+dc^*$. Thus, it is easily checked that
\begin{align}&\nn
ab^*+ba^*-ac^{-1}(cd^*+dc^*)(c^{-1})^*a^*+ac^{-1}cd^*q^*+qdc^*(c^{-1})^*a^*=0,
\end{align}
and \eqref{C50} takes the form
\begin{align}\nn
\clm_{22}\big(N, \ov{\zeta}, {\zeta}\big)^{-1}=&-q(N,\zeta)d(N,\zeta)c(N,\zeta)^*\clm_{22}\big(N,  {\zeta}, \ov{\zeta}\big)^{-1}
\\ & \label{C51} \times
c(N,\zeta)d(N,\zeta)^*q(N,\zeta)^* \qquad (\zeta \not=\ov{\zeta}, \,\, \zeta \not=\pm 2\I).
\end{align}
In view of the uniform limits \eqref{C10'} and \eqref{C44} and equality \eqref{C51}, we derive  (uniformly
for the values $\zeta$ on any compact $\clc \subset \BC_+ \backslash \{2\I\}$) the relation
\begin{align}& \label{C52}
\lim_{k\to \infty}\|q(k,\zeta)d(k,\zeta)c(k,\zeta)^*\|=0.
\end{align}

According to \eqref{C17} and Remark \ref{RkNond}, the inequality 
$$\Re\big(c(N,\xi)^{-1}d(N,\xi)\big)>0$$ 
is valid for $\xi \in \BC_+$.
Hence, we have
\begin{align}& & \label{C53}
2 \Re\big(c(N,\ze)c(N,\xi)^{-1}\clm_{22}\big(N,  {\xi}, \ov{\ze}\big)\big)>\clm_{22}\big(N,  {\ze}, \ov{\ze}\big) \quad (\xi \in \BC_+).
\end{align}
Using \eqref{C53} (and \eqref{C10+}) one easily proves (by negation) that for  any compact $\clc\subset \BC_+$ 
there are such $\wh N$ and $M >0$ that for each pair $\xi,\ze \in \clc$ and $N>\wh N$ the inequality
\begin{align}& \label{C54}
\|c(N,\xi)c(N,\ze)^{-1}\|<M
\end{align}
holds. It is immediate from \eqref{C52} and \eqref{C54} that we have a uniform limit
\begin{align}& \label{C55}
\lim_{k\to \infty}\|q(k,\zeta)d(k,\zeta)c(k,\xi)^*\|=0 \quad (\zeta, \xi \in \BC_+ \backslash \{2\I\}).
\end{align}

Taking into account \eqref{LF1}, \eqref{C40}, and uniform limits \eqref{C42} and \eqref{C43},  we see that
\begin{align}& \label{C46}
\lim_{k\to \infty}a(k,\zeta)c(k,\zeta)^{-1}=\lim_{k\to \infty}b(k,\zeta)d(k,\zeta)^{-1}=-\I \om(\zeta)
\end{align}
uniformly for the values $\zeta$ on any compact $\clc \subset \BC_+ \backslash \{2\I\}$.

Relations \eqref{C48}, \eqref{C55} and \eqref{C46} imply the equality
\begin{align}\nn
\lim_{k\to \infty}\clm_{12}(k,\zeta, \ov{\xi})&= a(N,\zeta)d(N,\xi)^*+ b(N,\zeta) c(N,\xi)^*
\\ & \label{C56}
=(-\I/2\pi) \om(\ze)\big(G_{\tau}(\xi)^*G_{\tau}(\zeta)\big)^{-1} \quad
(\zeta, \xi \in \BC_+ \backslash \{2\I\}),
\end{align}
which holds uniformly
on any compact $\clc \subset \BC_+ \backslash \{2\I\}$.
In this way, the reductions of \eqref{C30} for the blocks $\clm_{12}(k,\zeta,\ov{\xi})$ and $\clm_{21}(k,\zeta,\ov{\xi})=\clm_{12}(k,\xi,\ov{\zeta})^*$
of $\clm(k,\zeta,\ov{\xi})$ are proved.

Finally, using again \eqref{C48}, \eqref{C55} and \eqref{C46}, we see that the asymptotic equality
\begin{align}\nn
\lim_{k\to \infty}\clm_{11}(k,\zeta, \ov{\xi})&= a(N,\zeta)b(N,\xi)^*+ b(N,\zeta) a(N,\xi)^*
\\ & \label{C57}
=(1/2\pi) \om(\ze)\big(G_{\tau}(\xi)^*G_{\tau}(\zeta)\big)^{-1}\om(\xi)^* 
\end{align}
holds uniformly for $\zeta, \, \xi$  on any compact $\clc \subset \BC_+ \backslash \{2\I\}$).
Formulas \eqref{C10+}, \eqref{C56} and \eqref{C57} prove the theorem.
\end{proof}
\begin{Rk} \label{RkGen} We note that the asymptotic formula \eqref{C10+} for $\clm_{22}(k,\zeta, \ov{\xi})$ holds
uniformly for $\zeta, \xi$ belonging to the compacts $\clc$ from $\BC_+$ $($and not to $\clc \subset \BC_+ \backslash \{2\I\}).$
In the same way, uniform limits for $\clm_{11}(k,\zeta, \ov{\xi})$ and $\clm_{12}(k,\zeta, \ov{\xi})$ can be proved
as well. Moreover, the mentioned above proofs of \eqref{C30} admit generalizations for a wide class of interpolation 
problems $($see the abstract interpolation in \cite{SaLInt} and the proof of \eqref{C10} for  interpolation problems
in \cite{ALS87}$)$.
\end{Rk}
\begin{Rk}\label{RkPart} We note also that formula \eqref{C44} in the proof of Theorem \ref{TmM} is
of independent interest.
\end{Rk}

\section{Explicit formulas}\label{Expl}
\setcounter{equation}{0}

Let us consider the case, where the matrices $C_k$ are obtained explicitly,
namely, the case of the GBDT transformations of the trivial discrete Dirac system, where $C_k \equiv I_{2p}$.
Each GBDT is determined by some $n\in \BN$ and an {\it admissible} triple consisting of $n\times n$ matrices
$\cla$ and $\cls_0$ and of $n \times 2p$ matrix $\Pi_0$. The admissible triples are defined here as follows.
\begin{Dn} The triple $\{\cla,\cls_0,\Pi_0\}$ is called admissible if 
\begin{align}& \label{E1}
\cla\cls_0-\cls_0 \cla^*=\I \Pi_0 j \Pi_0^*, \quad \det \cla\not=0, \quad \cls_0>0.
\end{align}
\end{Dn}
According to \cite[Propositions 2.4, 3.1]{FKRS08}, each admissible triple determines 
{\it potential} $\{C_k\}$ $(k\geq 0)$ satisfying \eqref{I2} in the following way:
\begin{align} & \label{E2}
\Pi_{k+1}=\Pi_k+\I \cla^{-1}\Pi_k j \quad (k\geq 0),
\\ &  \label{E3}
S_{k+1}=S_k+\cla^{-1}S_k (\cla^*)^{-1}+\cla^{-1}\Pi_k
\Pi_k^*(\cla^*)^{-1} \quad (k\geq 0),
\\ &  \label{E4}
C_k:=I_{2p}+\Pi_k^*S_k^{-1}\Pi_k-\Pi_{k+1}^*S_{k+1}^{-1}\Pi_{k+1} .
\end{align}
We note that the matrices $\Pi_k$ in \eqref{E1}--\eqref{E4} have nothing to do with $\Pi(k)$ in Section \ref{CD}.
The Weyl function of Dirac system determined by the admissible triple $\{\cla,\cls_0,\Pi_0\}$ is given
by the formulas \cite[(4.7), (4.28)]{FKRS08}, which (in our notations) take the form 
\begin{align}& \label{E5}
\vp(\la)=-\I \big(I_p-\phi(\la)\big)\big(I_p+\phi(\la)\big)^{-1}, 
\\  & \label{E6}
 \phi(\la)=-\I \vt_1^*\cls_0^{-1}\big(\cla^{\times}-\la I_n\big)^{-1}\vt_2, \quad \cla^{\times}=\cla+\I \vt_2 \vt_2^*\cls_0^{-1},
\end{align}
where $\vt_i$ are $p\times p$ blocks of $\Pi_0$, that is, $\Pi_0=:\begin{bmatrix}\vt_1 & \vt_2 \end{bmatrix}$.
\begin{Tm}\label{ExplicitTm} Let Dirac system \eqref{I1}, \eqref{I2} be determined by some
admissible triple $\{\cla,\cls_0,\Pi_0\}$. Then, the set $\{\nu\}\cup \{s_{-k}\}$ $(0 \leq k<\infty)$ $($or, equivalently,
the $p\times p$ matrix $\nu=\nu^*$ and the semi-infinite Toeplitz matrix$)$ corresponding to this Dirac system
via \eqref{P4} is given explicitly by the formulas:
\begin{align}& \label{EM}
\nu=\Im\big(I_p+2\I \vt_1^*\cls_0^{-1}(\wt \cla+\I I_n)^{-1}\vt_2\big), \quad \wt \cla:=\cla+\I \vt_2(\vt_2-\vt_1)^*\cls_0^{-1};
\\ & \label{EM'}
s_0=2\Re\big(I_p+2\I \vt_1^*\cls_0^{-1}(\wt \cla+\I I_n)^{-1}\vt_2\big), 
\\ &  \label{EM+}
s_{-k}=2\I \vt_1^*\cls_0^{-1}(\wt \cla+\I I_n)^{-1}\clu^{k-1}(\clu-I_n)\vt_2 \quad (k\geq 1),
\end{align}
where $\Im$ and $\Re$ denote imaginary and real parts of the matrices, $\wt \cla+\I I_n$ is invertible and
\begin{align}& \label{E13}
\clu=(\wt \cla -\I I_n)(\wt \cla +\I I_n)^{-1}.
\end{align}
\end{Tm}
\begin{proof}.
Using \eqref{E5} and \eqref{E6}, we will write down $\vp(\la)$ in a more simple way.
Indeed, it is well-known from  system theory (and is easily checked directly) that
\begin{align}& \label{E7}
\big(I_{2p}-\I \vt_1^*\cls_0^{-1}\big(\cla^{\times}-\la I_n\big)^{-1}\vt_2\big)^{-1}=I_{2p}+\I \vt_1^*\cls_0^{-1}\big(\wt \cla-\la I_n\big)^{-1}\vt_2,
\end{align}
where $ \wt \cla=\cla^{\times}-\I \vt_2\vt_1^*\cls_0^{-1}$. Clearly this definition of $\wt \cla$ coincides with the one in formula
\eqref{EM}.
Moreover, the definition of $\wt \cla$
yields
\begin{align}&\nn
\I^2 \vt_1^*\cls_0^{-1}\big(\cla^{\times}-\la I_n\big)^{-1}\vt_2\vt_1^*\cls_0^{-1}\big(\wt \cla-\la I_n\big)^{-1}\vt_2
\\ &\nn
 =
\I \vt_1^*\cls_0^{-1}\big(\cla^{\times}-\la I_n\big)^{-1}(\cla^{\times}-\wt \cla)\big(\wt \cla-\la I_n\big)^{-1}\vt_2
\\ &  \label{E8}
=\I \vt_1^*\cls_0^{-1}\big(\wt \cla-\la I_n\big)^{-1}\vt_2-\I \vt_1^*\cls_0^{-1}\big(\cla^{\times}-\la I_n\big)^{-1}\vt_2.
\end{align}
From \eqref{E5}--\eqref{E8} we derive
\begin{align}& \label{E9}
\vp(\la)=-\I \big(I_p+2\I \vt_1^*\cls_0^{-1}\big(\wt \cla-\la I_n\big)^{-1}\vt_2\big),
\end{align}
and the representation
\begin{align} \label{E10}
\I \vp\Big(\I \frac{z+1}{z-1}\Big)&= I_p+2\I \vt_1^*\cls_0^{-1}\left(\wt \cla-\I \frac{z+1}{z-1} I_n\right)^{-1}\vt_2
\end{align}
follows. Since 
\begin{align}& \label{E11}
\Pi_0=\begin{bmatrix}\vt_1 & \vt_2 \end{bmatrix}, \quad  \cla^{\times}=\cla+\I \vt_2 \vt_2^*\cls_0^{-1}, \quad 
\wt \cla=\cla^{\times}-\I \vt_2\vt_1^*\cls_0^{-1},
\end{align}
the matrix identity in \eqref{E1} can be rewritten in the form
\begin{align}& \label{E12}
\wt \cla \cls_0 - \cls_0 \wt \cla^*=\I(\vt_1-\vt_2)(\vt_1-\vt_2)^*.
\end{align}
Taking into account that $\cls_0>0$ and \eqref{E12} is valid, we see
that $\s(\wt \cla)\subset \ov{\BC_+}$, where $\s$ means spectrum. In particular, $\det(\cla +\I I_n)\not=0$, and the
matrix $\cla +\I I_n$ is, indeed, invertible.

From \eqref{E13} and \eqref{E10} we derive that
\begin{align}\nn
\I \vp\Big(\I \frac{z+1}{z-1}\Big)&= I_p+2\I (z-1)\vt_1^*\cls_0^{-1}\big((z-1)\wt \cla-\I ({z+1}) I_n\big)^{-1}\vt_2
\\   \label{E14} &
=I_p-2\I (z-1)\vt_1^*\cls_0^{-1}(\wt \cla+\I I_n)^{-1}(I_n-z\clu)^{-1}\vt_2. 
\end{align}
Using equality $-(z-1)(I_n-z\clu)^{-1}=I_n+z(I_n-z\clu)^{-1}(\clu-I_n)$ and \eqref{E14}, we obtain a so called {\it realisation}
of  $\I \vp\Big(\I \frac{z+1}{z-1}\Big)$:
\begin{align} \nn
\I \vp\Big(\I \frac{z+1}{z-1}\Big)=&I_p+2\I \vt_1^*\cls_0^{-1}(\wt \cla+\I I_n)^{-1}\vt_2
+2\I z\vt_1^*\cls_0^{-1}(\wt \cla+\I I_n)^{-1}(I_n-z\clu)^{-1}
\\ & \label{E15}
\times(\clu-I_n)\vt_2.
\end{align}
Finally, relations \eqref{P4}, \eqref{P7-}, and \eqref{E15}  yield \eqref{EM}--\eqref{EM+}.

\end{proof}
\begin{Rk}\label{RkSmSp} In view of \eqref{E13},  $\clu$ is invertible in the case $\I\not\in \s(\wt \cla)$.
In this case, the Toeplitz matrices considered in Theorem \ref{ExplicitTm} $($and given by the relations
\eqref{EM'}, \eqref{EM+} and $s_k=s_{-k}^*)$ are block diagonal plus block
semiseparable matrices.
\end{Rk}
We note that semiseparable matrices have been studied in various papers (see, e.g., \cite{EG, GKK, Van}
and the references therein).

\noindent{\bf Acknowledgments.}
 {This research    was supported by the
Austrian Science Fund (FWF) under Grant  No. P29177.}

\begin{flushright}
A.L. Sakhnovich,\\
Faculty of Mathematics,
University
of
Vienna, \\
Oskar-Morgenstern-Platz 1, A-1090 Vienna,
Austria, \\
e-mail: {\tt oleksandr.sakhnovych@univie.ac.at}

\end{flushright}

\end{document}